\documentclass[11pt]{article}
\usepackage{amsmath}
\usepackage{mathrsfs}
\usepackage{amsfonts}

\usepackage{amsfonts,amsmath,amssymb}
\usepackage{amssymb,amsfonts,amsmath,
latexsym,epsfig,cite,psfrag,eepic,colordvi}
\usepackage{amscd,graphics}
\newcommand{\qed}{\hfill $\Box $}
\newcommand{\pf}{\noindent {\bf Proof.} }

\newtheorem{theorem}{Theorem}[section]

\begin{document}

\begin{center}
{\Large\bf A Note on $\{k,n-k\}$-Factors of Regular Graphs}
\end{center}

\begin{center}
Hongliang Lu$^{1}$ and David G. L. Wang$^{2}$\\[6pt]

$^{1}$Department of Mathematics\\
$^{1}$Xi'an Jiaotong University, Xi'an 710049, P. R. China\\
{\tt $^{1}$luhongliang215@sina.com}

$^{2}$Beijing International Center for Mathematical Research\\
$^{2}$Peking University, Beijing 100871, P. R. China\\
{\tt $^{2}$wgl@math.pku.edu.cn}
\end{center}

\begin{abstract}
Let $r$ be an odd integer,
and $k$ an even integer.
In this note,
we present $r$-regular graphs
which have no $\{k,r-k\}$-factors for all $1\le k\le {r\over2}-1$.
This gives a negative answer to a problem posed
by Akbari and Kano recently.
\end{abstract}


\noindent\textbf{2010 AMS Classification:} 05C75

\section{Introduction}

Let~$G$ be a simple graph,
and~$H$ a set of integers.
A spanning subgraph~$F$ of~$G$ is called
an \emph{$H$-factor} if the degree of each vertex
belongs to~$H$.

Let $r$, $k$ be positive integers
such that $1\le k\le {r\over2}$.
Recently, Akbari and Kano~\cite{AK} considered
whether a general $r$-regular graph
contains a $\{k,r-k\}$-factor.
It is known~\cite{Lovasz} that
the problem of determining whether~$G$ has an $\{a,b\}$-factor
is NP-complete for any positive integers $a$ and $b$ such that $b-a\ge 3$.

The case that both $r$ and $k$ are even was solved by
Peterson~\cite{Pet1891} early in 1891. In fact Peterson proved that
any $r$-regular graph has a $\{k\}$-factor. The proof can be found
in Akiyama and Kano's book~\cite{AK}. For the case that $r$ is odd
and $k$ is even, Akbari and Kano~\cite{AK11} showed the existence of
a $\{k,r-k\}$-factor of any $r$-regular graphs by induction on
edge-connectivity. In fact, they considered general graphs (allowing
multiple edges and loops) and obtained some refined results. They
posed the problem of solving the remaining cases.

In this note,
we deal with the case that $r$ is even and $k$ is odd.

\section{The main result}

Let $r$ be an even integer, and $k$ an odd integer,
such that $1\le k\le {r\over2}$.
In this section,
we consider the problem that whether every $r$-regular graph
contains a $\{k,r-k\}$-factor.

First, the case $k={r\over 2}$ can be completely solved
by the following result of Gallai~\cite{Gal50}.

\begin{theorem}[Gallai]\label{thm_Gallai}
Every $m$-edge-connected $r$-regular graph on
an even number of vertices has a $k$-factor
if $r$ is even, $k$ is odd, and
\[
{r\over m}\le k\le r\biggl(1-{1\over m}\biggr).
\]
\end{theorem}

\begin{theorem}
Let $r$ be an even integer such that ${r\over2}$ is odd.
Let $G$ be a connected $r$-regular graph.
Then $G$ has an ${r\over2}$-factor if and only if
the number of vertices of $G$ is even.
\end{theorem}

\pf The necessity is clear since ${r\over2}$ is odd.
Suppose that $G$ is an $r$-regular graph with an even number of vertices.
Note that any connected $r$-regular graph is
$2$-edge-connected if $r$ is even.
Taking $m=2$ in Gallai's Theorem~\ref{thm_Gallai},
we see that $G$ has an ${r\over2}$-factor. \qed

For the other cases, we will present $r$-regular graphs
without $\{k,r-k\}$-factors.
Here is our main result.

\begin{theorem}
Let $r$, $k$ be positive integers such that
$r$ is even, $k$ is odd, and
\[
1\le k\le {r\over2}-1.
\]
Then there exists an
$r$-regular graph without $\{k,r-k\}$-factors.
\end{theorem}

\pf Note that the case~${r\over2}$ is even and $k={r\over2}-1$ is
solved by the first author~\cite{Lu}.  We divide the remaining cases
into two parts depending on the parity of~${r\over2}$. For either of
them we shall present an $r$-regular graph without
$\{k,r-k\}$-factors.

Let~$H$ be the graph obtained by removing
an edge from the complete graph~$K_{r+1}$.
We will use $H$ as building blocks in the constructions.
Let $H_1,H_2,\ldots$ be pairwise disjoint copies of~$H$.
Let~$u_i$ and~$v_i$ be the unsaturated vertices of~$H_i$.

Suppose that ${r\over2}$ is odd. Then $k\le {r\over2}-2$.
Let $G_1$ be the graph obtained by linking
a vertex $v$ to all unsaturated vertices
of $H_1$, $H_2$, $\ldots$, $H_{r\over2}$.
Then $G_1$ is an $r$-regular graph with $r(r+1)/2+1$ vertices.
%\begin{figure}[h,t]
%\begin{center}
%\includegraphics[width=6cm]{2.eps}
%\end{center}
%\vspace{-10mm}
% \caption{ }
%\end{figure}
Suppose that $G_1$ contains a $\{k,r-k\}$-factor $F_1$.
By parity, it is easy to deduce that
$v$ links to exactly one vertex of each $H_i$ in $F_1$.
It follows that
\[
\deg_{F_1}(v)={r\over2}\not\in\{k,r-k\},
\]
a contradiction. So $G_1$ is a desired graph.

Now we consider the case that ${r\over2}$ is even.
By the work in~\cite{Lu},
we can suppose that $k\le {r\over2}-3$.
In this case,
let $G_2$ be the graph consisting of two vertices $u$, $v$,
and the copies $H_1$, $H_2$, $\ldots$, $H_r$ such that
$u$ is linked to all unsaturated vertices
of the copies $H_1$, $H_2$, $\ldots$, $H_{{r\over2}-1}$,
and to the vertices $u_{r-1}$ and $u_r$;
while~$v$ is linked to all unsaturated vertices of
$H_{{r\over2}}$, $H_{{r\over2}+1}$, $\ldots$, $H_{r-2}$,
and to the vertices $v_{r-1}$ and $v_r$.
Then $G_2$ is an $r$-regular graph with $r(r+1)+2$ vertices.
%\begin{figure}[h,t]
%\begin{center}
%\includegraphics[width=8cm]{1.eps}
%\end{center}
%%\vspace{-10mm}
%\caption{a}
%\end{figure}
Suppose that $G_2$ contains a $\{k,r-k\}$-factor $F_2$.
By parity, the vertex $u$ links to exactly one vertex of each $H_i$
($1\le i\le {r\over2}-1$). It follows that
\[
\deg_{F_2}(u)\in\biggl\{{r\over2}-1,\ {r\over2},\ {r\over2}+1\biggr\},
\]
which has empty intersection with the set $\{k,r-k\}$. Hence $G_2$
has no $\{k,r-k\}$-factors. This completes the proof.\qed

\end{document}